\documentclass[a4paper,12pt]{article}
\catcode`\"=13
\def"#1{\v#1}
\usepackage{graphicx,amsmath,amssymb,euscript,amsfonts}
\usepackage{textcomp}
\graphicspath{{images/}}
\usepackage[english]{babel}

\newcommand{\qed}{\hfill\hbox{\rule{3pt}{6pt}}}
\newcommand{\proof}{{\sc Proof. }}

\newtheorem{theorem}{Theorem}[section]
\newtheorem{lemma}[theorem]{Lemma}

\newtheorem{remark}[theorem]{Remark}
\title{On a problem of \\
Bermond and Bollob\'as}

\author{
         Slobodan Filipovski\footnote{Supported in part by the Slovenian Research Agency (research program P1-0285 and Young Researchers Grant).} \\
         University of Primorska, Koper, Slovenia \\
         \texttt{slobodan.filipovski@famnit.upr.si}
         \and
         Robert Jajcay\thanks{Supported by VEGA 1/0474/15, VEGA 1/0596/17, APVV-15-0220, Slovenian Research Agency J1-6720, and by NSFC 11371307.}\\
        Comenius University, Bratislava, Slovakia \\
         \texttt{robert.jajcay@fmph.uniba.sk}
 }

\begin{document}
 \maketitle

\begin{abstract}
Let $n(k, d)$ be the order of the largest undirected graphs of
maximum degree $k$ and diameter $d$, and let $M(k,d)$ be the corresponding
Moore bound.
In this paper, we give a positive answer to the question of Bermond and
Bollob\'as concerning the Degree/Diameter Problem:
Given a positive integer $c>0$, does there exist a pair $k$ and $d$, such that $n(k, d)\leq M(k,d)-c?$

 \end{abstract}

\section{Introduction} \label{sec:intro}
\medskip

\noindent
We call a $k$-regular graph $\Gamma$ of diameter $d$ a \emph{$(k,d)$-graph}.
Let $n(k,d)$ denote the largest order of any undirected graph of
maximum degree $k$ and diameter $d$. It is easy to show that the
order $|V(\Gamma)|$ of any graph $\Gamma$ of maximum degree $k$ and diameter
$d$ and therefore also the parameter $n(k,d)$ satisfy the following
inequality:
$$|V(\Gamma)| \leq n(k,d)\leq M(k,d)=1+k+k(k-1)+k(k-1)^{2}+...+k(k-1)^{d-1}.$$
The above number $M(k,d)$ is called the \emph{Moore bound}. A graph whose order is equal to the Moore bound is called a \emph{Moore graph}; such a graph is necessarily regular of degree $k$. Moore graphs are proved to be very rare. They are the complete graphs on $k+1$ vertices; the cycles on $2d+1$ vertices; and for diameter $2$, the Petersen graph, the Hoffman-Singleton graph and possibly a graph of degree $k=57$.
For $k>2$ and $d>2$, there are no Moore graphs \cite{K}.

The main problem addressed in this paper is closely related to the
well-known extremal graph theory problem called
the \emph{Degree/Diameter Problem}:
\begin{center}
\emph{Given natural numbers $k$ and $d$, find the largest possible number of vertices $n(k,d)$ in a graph of maximum degree $k$ and diameter at most
$ d$.}
\end{center}
The difference between the Moore bound $M(k,d)$ and the order of a
specific graph $\Gamma$ of maximum degree $k$ and diameter $d$ is called the \emph{defect} of $\Gamma$,
and is denoted by $\delta(\Gamma)$. Thus, if $\Gamma$ is a largest
graph of maximum degree $k$ and diameter $d$, then $n(k,d)=M(k,d)-\delta(\Gamma)$. It needs to be
noted that very little
is known about the exact relation between the Moore bounds $M(k,d)$
and the corresponding extremal orders $n(k,d)$. While a considerable
gap exists between
the orders of the largest known/constructed graphs of maximum
degree $k$ and diameter $d$ and the corresponding Moore bounds, no
substantially better bounds are known. In particular, it is not even known whether the two parameters are of the same order of magnitude (with computational evidence strongly suggesting that they are not).
Thus, the below stated long open {\em question of Bermond and
Bollob\' as} \cite{B} can be
viewed as the natural first attempt at shedding light on the nature of the
relation between $M(k,d)$ and $n(k,d)$. In this paper,
we answer in positive the following:
\begin{center}
\emph{Is it true that for each positive integer $c$ there exist $k$ and $d$ such that the order of the largest graph of maximum degree $k$ and diameter $d$ is at most $M(k,d)-c?$}
\end{center}

Moore graphs of degree $k$ and diameter $d$ are well-known to be the
only $(k,d)$-graphs of girth $2d+1$; the girth of any other (i.e., non-Moore) $(k,d)$-graph is strictly smaller than $2d+1$. It will also prove
useful to note that even though it is not known whether the extremal
graphs of diameter $d$, maximal degree $k$, and of the maximal order $n(k,d)$, are necessarily $k$-regular, the graphs $\Gamma$ conatining vertices of
degree smaller than $k$ must satisfy the following stricter upper bound:
\begin{equation}\label{trivial}
|V(\Gamma)| \leq M(k,d) - 1 - (k-1) - \ldots - (k-1)^{d-1} =
M(k,d) - \frac{(k-1)^d-1}{k-2}  .
\end{equation}
The \emph{bipartite Moore bound} is the maximum number $B(k,d)$ of vertices in a bipartite graph of maximum degree $k$ and diameter at most $d$. This bound
is due to Biggs \cite{L}:
\begin{equation}\label{bipart}
B(2,d)=2d, \mbox{  and  }  B(k,d)=\dfrac{2(k-1)^{d}-2}{k-2},
\mbox{ if } k>2,
\end{equation}
and is smaller than the Moore bound by $(k-1)^d$ (i.e, $M(k,d)-B(k,d) =
(k-1)^d$, for $ k \geq 3$).
Bipartite $(k,d)$-graphs of order $B(k,d)$ are called \emph{bipartite Moore graphs}.
The bipartite Moore bound represents not only an upper bound on the number of vertices of a bipartite graph of maximum degree $k$ and diameter $d$,
but it is also a lower bound on the number of vertices of a regular graph
$\Gamma$ of degree $k$ and girth $g=2d$. A $(k,g)$\emph{-cage} is a
smallest $k$-regular graph of girth $g$, and if a $(k,2d)$-cage is of order $B(k,d)$,
it is a bipartite Moore graph \cite{L}.

For degrees $1$ or $2$, bipartite Moore graphs consist of $K_{2}$ and the $2d$-cycles, respectively. When $k\geq3$, the possibility of the existence of bipartite Moore graphs was settled by Feit and Higman \cite{F}  in $1964$ and, independently, by Singleton \cite{S} in $1966$. They proved that such graphs exist only for the diameters $2, 3, 4$ or $6$.

The question of Bermond and Bollob\' as has already been answered in positive for
the more specialized families of vertex-transitive and Cayley graphs \cite{A}.
A graph $\Gamma$ is \emph{vertex-transitive} if its automorphism group acts
transitively on its set of vertices $V(\Gamma)$, i.e., for every pair of vertices $u,v
\in V(\Gamma)$ there exists an automorphism $\varphi$
of $\Gamma$ mapping $u$ to $v$, $\varphi(u)=v$. A vertex-transitive graph
$\Gamma$ is said to be \emph{Cayley} if its automorphism group contains a subgroup $G$ acting regularly on $V(\Gamma)$, i.e., having the property that
for every pair of vertices $u,v \in V(\Gamma)$ there exists exactly one
automorphism $\varphi \in G$ mapping $u$ to $v$, $\varphi(u)=v$.
Due to their high level of symmetry, all vertex-transitive graphs are regular.
If we let $vt(k,d)$ denote the largest order of a vertex-transitive $(k,d)$-graph, and $C(k,d)$ denote the largest order of a Cayley graph of degree $k$ and diameter $d$, then
\[ C(k,d)\leq vt(k,d) \leq n(k,d) . \]
Exoo et al. proved in \cite{A} that for any fixed $k\geq3$ and $c\geq2$ there exists a set ${\mathcal S}$ of natural numbers of positive density such that $vt(k,d)\leq M(k,d)-c,$ for all $d\in {\mathcal S}.$ The same holds for Cayley graphs as well.
We list the result for the sake of completeness.

\begin{lemma}[\cite{A}] Let $k\geq3 $ and $c\geq2$. Let $r$ be an odd integer, and let $p$ be a prime such that $p>2k(k-1)^{(r-1)/2}>8k(k-1)^{2}c^{2}.$ If $2d+1=rp,$ then any vertex-transitive $(k,d)$-graph has defect greater than $c$.
\end{lemma}
In this paper, we do not seek to obtain density results. Rather, we prove that for any $ c > 0 $,
there exist infinitely many pairs $(k,d)$ for which $ n(k,d) <  M(k,d)-c $.
Our approach is a combination of applications of previously known methods and of some new uses of spectral analysis.



\section{Spectral analysis}
For many of our arguments, we rely on the techniques of
spectral analysis applied to graphs extremal with respect to
the Degree/Diameter Problem.
We begin with a brief review of the basic facts as listed in \cite{K}.

Let $\Gamma$ be a connected $(k,d)$-graph of order $n$ and defect $\delta>0.$
For each integer $i$ in the range $0\leq i\leq d$, we define the $n\times n$ \emph{$i$-distance matrices } ${\bf A}_{i}={\bf A}_{i}(\Gamma)$ as follows. The rows and columns of ${\bf A}_{i}$ correspond to the vertices of $\Gamma$, and the entry in position $(\alpha, \beta)$ is $1$ if the distance between the vertices $\alpha$ and $\beta$ is $i$, and zero otherwise. Clearly
${\bf A}_{0}={\bf I}$ and ${\bf A}_{1}={\bf A}$, the usual adjacency matrix of $\Gamma$. If ${\bf J}_{n}$ is the all-ones matrix and $d$ is the diameter of a connected graph $\Gamma$, then
$$ \sum_{i=0}^{d}{\bf A}_{i}={\bf J}_{n}. $$
Given a finite $k$-regular graph $\Gamma$,
we define the polynomials $G_{k,i}(x)$ for all $x \in \mathbb{R}$
recursively as follows:
\begin{equation}\label{eq}
\left\{
\begin{array}{lc}G_{k,0}(x)=1  \\
G_{k,1}(x)=x+1 \\
G_{k,i+1}(x)=xG_{k,i}(x)-(k-1)G_{k, i-1}(x)\mbox{ for } i\geq1.
\end{array}
\right.
\end{equation}
We note that the entry $(G_{k,i}({\bf A}))_{\alpha, \beta}$ counts the number of paths of length at most $i$ joining the vertices $\alpha$ and $\beta$ in $\Gamma$.
Regular graphs with defect $\delta$ and order $n$ satisfy the matrix equation
$$G_{k,d}({\bf A})={\bf J}_{n}+{\bf B},$$
where ${\bf B}$ is a non-negative integer matrix with the row and column sums equal to $\delta.$ The matrix ${\bf B}$ is called the \emph{defect matrix}, (see \cite{K}).

Next, we follow the line of argument that originally appeared in \cite{Big}.
Since $\Gamma$ is regular and connected, the all-ones matrix ${\bf J}_{n}$ is a polynomial of
${\bf A}$, say, $ J_n({\bf A}) $. From now on, we adopt the convention that matrices will be denoted by upper-case bold-face characters while their corresponding polynomials will be denoted by the same character but not bold-faced. Thus,  ${\bf B}=B({\bf A})=G_{k,d}({\bf A})-J_n({\bf A})$, and
$ {\bf J}_n = J_n({\bf A}) = G_{k,d}({\bf A}) - B({\bf A}) $.
It follows that if $\lambda$ is an eigenvalue of ${\bf A}$, then $G_{k,d}(\lambda)-B(\lambda) = J_n(\lambda) $ is an eigenvalue of ${\bf J}_{n}.$ Substituting the value $k$ for $\lambda$ yields the eigenvalue $n$ of ${\bf J}_{n}$, $G_{k,d}(k)-B(k) = n $. An easy calculation yields
that $G_{k,d}(k)=M(k,d)$, and therefore $B(k)=M(k,d)-n=\delta$ is an eigenvalue of ${\bf B}$. Since each row and column of ${\bf B}$ sums up to $\delta$, every eigenvalue of ${\bf B}$ has value at most $\delta.$
If $\lambda\neq k$ is another eigenvalue of ${\bf A}$, then $G_{k,d}(\lambda)-B(\lambda)$ must be the zero eigenvalue of ${\bf J}_{n}$.
Therefore, $G_{k,d}(\lambda)-B(\lambda) = 0 $, and since $|B(\lambda)|\leq \delta$, we obtain
$|G_{k,d}(\lambda)|\leq \delta.$
Thus, the value $|G_{k,d}(\lambda)|$ is a lower bound for the defect
$ \delta(\Gamma)$. In summary, if $\Gamma$ is a graph of diameter $d$, degree $k$, and order $M(k,d)-\delta$, then every eigenvalue $\lambda \neq k$ of $\Gamma$ satisfies
\begin{equation} \label{1}
|G_{k,d}(\lambda)|\leq \delta.
\end{equation}

Since ${\bf A}$ is symmetric, all eigenvalues of ${\bf A}$ are real. Let $\lambda_{0} \geq \lambda_{1} \geq \ldots \geq \lambda_{n-1}$ be the eigenvalues of $\Gamma$ and let $\lambda$ be the eigenvalue with the  second largest absolute value. It is well known from \emph{Perron-Frobenius theory} (e.g., \cite{Lov}), that $\lambda_{0}=k$ and $k$ is of multiplicity one if and only if $\Gamma$ is connected. Moreover, if $\Gamma$ is non-bipartite, then $\lambda_{n-1}> -k$. Therefore, if $\Gamma$ is a connected and non-bipartite graph, $\lambda=\max\{\lambda_1, |\lambda_{n-1}|\}.$ Studying the second largest eigenvalue of a given graph is related to the existence of the \emph{Ramanujan graphs}. We call a $k$-regular  graph a \emph{Ramanujan graph} if its second largest eigenvalue (in absolute value) is at most $2\sqrt{k-1}$
\cite{ramanujan}.

In this paper we use a result obtained by Alon and Boppana in \cite{AB}, which provides a lower bound on $\lambda$  for graphs of very large order. Let $X_{n,k}$ be a $k$-regular graph on $n$ vertices and let $\lambda(X_{n,k})$ denote the eigenvalue of the second largest absolute
value.
Our proof is based on the following restatement of this result.

\begin{theorem}[\cite{ramanujan}] \label{important} $ \liminf_{n\rightarrow \infty} \lambda(X_{n, k})\geq 2 \sqrt{k-1}.$
\end{theorem}

\newpage

\section{Main result}

The following theorem is the main result of this paper.
\begin{theorem} Let $c>0$ be a fixed integer and let $k \geq 3 $ be larger
than $c$. Then there exists at least one large enough even $d > k$ such that any $(k,d)$-graph $\Gamma$ has defect greater than $c$, i.e.,   $n(k,d)\leq M(k,d)-c.$
\end{theorem}

\noindent
\proof
Let $ c \geq 1 $ and $ k \geq 3 $, $ k > c $. In order to prove the
theorem, we divide the class of graphs of degree
not exceeding $k$ into disjoint subclasses which we then treat
separately, choosing subsequently larger and larger diameters
$d$ satisfying the desired claim until we find a diameter $d$ that
works for the entire class.

As argued in the introduction, Moore and bipartite Moore graphs
exist only for very limited diameters $d$. Therefore, assuming that
$d > 6 $, yields the non-existence of Moore or bipartite Moore $(k,d)$-graphs (and hence the non-existence of $(k,d)$-graphs of girth
$2d+1$). In addition, due to (\ref{trivial}), taking $ d > \log_{k-1}(c(k-2)+1)$
makes the order of any {\em non-regular} graph of maximum degree $k$ and diameter $d$ smaller than $ M(k,d)-c $. Similarly, due to (\ref{bipart}),
taking $ d > \log_{k-1}(c) $ makes the order of any {\em bipartite} $(k,d)$-graph
smaller than $ M(k,d)-c $. Because of our assumptions on $c$ and
$k$, both $ \log_{k-1}(c(k-2)+1) $ and $\log_{k-1}(c)$ are smaller than
$6$, and thus, considering $ d >  6 $ yields the {\em non-existence} of
graphs of orders larger than $M(k,d)-c$ in the families of $(k,d)$-graphs
of girth $2d+1$, the non-regular graphs
of maximum degree $k$ and diameter $d$, and the
bipartite $(k,d)$-graphs.
Therefore, from now on,
we will assume that $d>6$ and $\Gamma$ is $k$-regular,
non-bipartite and of girth smaller than $2d+1$.

The proof of the non-existence of such graphs of orders larger
than $ M(k,d)-c$ for sufficiently large $d$'s splits into two cases: The case when $\Gamma$ is a non-bipartite $(k,d)$-graph of girth $2d$ or $2d-1$, and the case when $\Gamma$ is a non-bipartite $(k,d)$-graph of girth smaller than $2d-1$.


Let us begin with the case when the girth of $\Gamma$ is $2d$ or $2d-1$. In this case, it is easy to see that
the order of $\Gamma$ satisfies $|V(\Gamma)| \geq 1+k+k(k-1)+\ldots+k(k-1)^{d-2}$. This clearly implies that the order of $(k,d)$-graphs
$\Gamma$ increases to infinity when $d$ does.
Theorem \ref {important} yields that for a fixed $\delta>0,$ there exists $N>0$, such that for all graphs $\Gamma$ on $n>N$ vertices holds $\lambda(\Gamma) \geq 2 \sqrt{k-1}-\delta.$
Thus, for a fixed $\delta>0$ all $(k,d)$-graphs $\Gamma$ with arbitrary large diameter satisfy at least one of the conditions
$\lambda_{1}\in [2\sqrt{k-1}-\delta, k)$ or $\lambda_{n-1} \in (-k, -2\sqrt{k-1}+\delta]$
(recall that $\Gamma$ is assumed connected and non-bipartite, hence the half-open intervals).
Since $G_{k,d}(x)$ is a continuous function on the closed interval $[-k,k]$, it follows that $G_{k,d}(x)$ is uniformly continuous on $[-k,k]$.
It implies that for every $\epsilon>0$ there exists $\delta>0$,  such that for every $x, x_{0}\in [-k, k]$ holds
$$| x-x_{0}|<\delta \Rightarrow |G_{k,d}(x)-G_{k,d}(x_{0})|<\epsilon.$$
 Therefore, fixing $x_{0}=2\sqrt{k-1}$ and a small $\epsilon>0$, it implies the existence of $\delta>0$ such that for $x\in(2\sqrt{k-1}-\delta, 2\sqrt{k-1}]$ we have $G_{k,d}(x)>G_{k,d}(2\sqrt{k-1})-\epsilon.$ Similarly we obtain $G_{k,d}(x)>G_{k,d}(-2\sqrt{k-1})-\epsilon$ for $x\in[-2\sqrt{k-1}, -2\sqrt{k-1}+\delta).$ These conclusions lead to the fact that the values of $G_{k,d}(x)$ on $[-2\sqrt{k-1}, -2\sqrt{k-1}+\delta)\cup(2\sqrt{k-1}-\delta, 2\sqrt{k-1}]$ can be estimated based on the values $G_{k,d}(\pm2\sqrt{k-1})$.

 In order to take advantage of (\ref{1}),
we will derive explicit formula for $G_{k,d}(x)$.
Fixing the variable $x$ makes the last equation of (\ref{eq}) into a second order linear homogeneous recurrence equation for $G_{k,d}(x)$ with respect to the
parameter $d$ subject to the initial conditions $G_{k,0}(x)=1$
and $G_{k,1}(x)=x+1$. Since we only need to calculate the values of
$G_{k,d}(x)$ for $ x \in (-k, -2\sqrt{k-1}]  \cup [2\sqrt{k-1}, k)$, we only
need to consider the recurrence relation in the case when the roots
of the corresponding second degree polynomial equation
$ t^2 -xt +(k-1) = 0 $ are real, with a double-root when $ x = \pm 2 \sqrt{k-1} $.
Solving this recurrence equation for a fixed
$x \in (-k, -2\sqrt{k-1})  \cup (2\sqrt{k-1}, k)$, we obtain the explicit formula
 $$G_{k,d}(x)=\frac{x+2+\sqrt{x^{2}-4k+4}}{2\sqrt{x^{2}-4k+4}}\left(\frac{x+\sqrt{x^{2}-4k+4}}{2}\right)^{d}-$$
 $$-\frac{x+2-\sqrt{x^{2}-4k+4}}{2\sqrt{x^{2}-4k+4}}\left(\frac{x-\sqrt{x^{2}-4k+4}}{2}\right)^{d}.$$

In the case when $x=\pm2\sqrt{k-1}$, the second degree polynomial
equation has a double root and we obtain:
$$ G_{k,d}(2\sqrt{k-1})=(d+1)\sqrt{k-1}^{d}+d\sqrt{k-1}^{d-1},$$
$$ G_{k,d}(-2\sqrt{k-1})=(-1)^{d}((d+1)\sqrt{k-1}^{d}-d\sqrt{k-1}^{d-1}).$$

It is easy to see from (\ref{eq}) that the function $G_{k,d}(x)$ is a
polynomial of degree $d$ in $x$, and thus differentiable.
Calculating the derivative of $G_{k,d}(x)$, for
$x \in (-k, -2\sqrt{k-1})  \cup (2\sqrt{k-1}, k)$, we have
$$G_{k,d}^{'}(x)=\frac{d(x-\sqrt{x^{2}-4k+4})^{d}(x+2-\sqrt{x^{2}-4k+4})}{2^{d+1}(x^{2}-4k+4)}+\frac{(x-\sqrt{x^{2}-4k+4})^{d+1}}{2^{d+1}(x^{2}-4k+4)}+$$
$$+\frac{x(x-\sqrt{x^{2}-4k+4})^{d}(x+2-\sqrt{x^{2}-4k+4})}{2^{d+1}(x^{2}-4k+4)^{\frac{3}{2}}}+\frac{d(x+\sqrt{x^{2}-4k+4})^{d}(x+2+\sqrt{x^{2}-4k+4})}{2^{d+1}(x^{2}-4k+4)}+$$
$$+\frac{(x+\sqrt{x^{2}-4k+4})^{d+1}}{2^{d+1}(x^{2}-4k+4)}-\frac{x(x+\sqrt{x^{2}-4k+4})^{d}(x+2+\sqrt{x^{2}-4k+4})}{2^{d+1}(x^{2}-4k+4)^{\frac{3}{2}}}
.$$
The behavior of $G_{k,d}(x)$ differs in the intervals
$ (2\sqrt{k-1}, k) $ and $ (-k, -2\sqrt{k-1})$.
First, let us suppose that $\lambda_{1}\in (2\sqrt{k-1}, k)$. The
assumption $x\in(2\sqrt{k-1}, k)$ yields the inequalities
$x^{2}-4k+4>0, x>\sqrt{x^{2}-4k+4}$, and $d>k>x$.

We note that five of the terms of $G_{k,d}^{'}(x)$ are positive numbers. Clearly, if $x^{2}-4k+4\geq1$, then
$\frac{d(x+\sqrt{x^{2}-4k+4})^{d}(x+2+\sqrt{x^{2}-4k+4})}{2^{d+1}(x^{2}-4k+4)}>\frac{x(x+\sqrt{x^{2}-4k+4})^{d}(x+2+\sqrt{x^{2}-4k+4})}{2^{d+1}(x^{2}-4k+4)^{\frac{3}{2}}}$
and thus we easily see that $G_{k,d}^{'}(x)>0.$\\
Now, let us suppose that $0<x^{2}-4k+4<1.$
If $x^{2}-4k+4\rightarrow0$, that is, if $x\rightarrow2\sqrt{k-1}$, then
$$\lim_{x\rightarrow2\sqrt{k-1}}\frac{x(x-\sqrt{x^{2}-4k+4})^{d}(x+2-\sqrt{x^{2}-4k+4})}{2^{d+1}(x^{2}-4k+4)^{\frac{3}{2}}}=$$
$$=\lim_{x\rightarrow2\sqrt{k-1}}\frac{x(x+\sqrt{x^{2}-4k+4})^{d}(x+2+\sqrt{x^{2}-4k+4})}{2^{d+1}(x^{2}-4k+4)^{\frac{3}{2}}}.$$
Moreover, if $x\rightarrow2\sqrt{k-1}$, then the other four terms of $G_{k,d}^{'}(x)$ are positive and tend to infinity.
Thus, we get $G_{k,d}^{'}(x)>0$.\\
If $0<x^{2}-4k+4<1$ and $x^{2}-4k+4\nrightarrow0$, then we always can choose large enough $d$ such that $d>\frac{x}{\sqrt{x^{2}-4k+4}}.$
It implies the inequality $\frac{d(x+\sqrt{x^{2}-4k+4})^{d}(x+2+\sqrt{x^{2}-4k+4})}{2^{d+1}(x^{2}-4k+4)}>\frac{x(x+\sqrt{x^{2}-4k+4})^{d}(x+2+\sqrt{x^{2}-4k+4})}{2^{d+1}(x^{2}-4k+4)^{\frac{3}{2}}}$,
and hence $G_{k,d}^{'}(x)>0$.

 Using these inequalities we can deduce that the functions $G_{k,d}(x)$ and $G_{k,d}^{'}(x)$  are both positive on the interval $(2\sqrt{k-1}, k)$. Since
$G_{k,d}(x)$ is continuous on $[2\sqrt{k-1}, k]$ (being a polynomial), and differentiable on $(2\sqrt{k-1}, k)$,
$G_{k,d}(x)$ is increasing on $[2\sqrt{k-1}, k]$. This provides us
with a lower bound on the defect of $\Gamma$ in this case:
$$\delta(\Gamma)\geq |G_{k,d}(\lambda_{1})|=G_{k,d}(\lambda_{1})\geq G_{k,d}(2\sqrt{k-1})=(d+1)\sqrt{k-1}^{d}+d\sqrt{k-1}^{d-1}>c.$$
Next, let us suppose that $\lambda_{n-1}\in (-k, -2\sqrt{k-1})$. If $x\in (-k, -2\sqrt{k-1})$, then $0>x+2+\sqrt{x^{2}-4k+4}>x+2-\sqrt{x^{2}-4k+4}$,
and since $d$ is assumed to be an {\em even} number,
$(x-\sqrt{x^{2}-4k+4})^{d}>(x+\sqrt{x^{2}-4k+4})^{d}>0$.
Based on these inequalities, we observe that the function $G_{k,d}(x)$ is positive on $(-k, -2\sqrt{k-1}).$
Moreover, using the above inequalities again, we can deduce that $G_{k,d}^{'}(x)$ is negative on $(-k, -2\sqrt{k-1}).$
Thus, if $d$ is an even number, then $G_{k,d}(x)$ is decreasing on $[-k, -2\sqrt{k-1}].$ This implies the bound
\begin{eqnarray*}
\delta(\Gamma)\geq |G_{k,d}(\lambda_{n-1})|=
G_{k,d}(\lambda_{n-1})\geq G_{k,d}(-2\sqrt{k-1})= \\
= (-1)^{d}((d+1)\sqrt{k-1}^{d}-d\sqrt{k-1}^{d-1})>c.
\end{eqnarray*}
The two above inequalities yield that any non-bipartite $(k,d)$-graph
$ \Gamma $ of diameter $ d > D $ and girth $ 2d $ or $ 2d-1 $ is of
defect large than $c$.



Finally, let us assume that the girth $g(\Gamma)$ of $\Gamma$ is at most $2d-2$, i.e., $3\leq g(\Gamma)\leq 2d-2$.
The main idea of this part of the proof is the observation that a small girth forces a loss
of an entire branch of the potential Moore tree of $(k,d)$-graph with the number of `lost' vertices
being roughly of the order of magnitude of $(k-1)^{d-g/2}$.
Assume that $b \in V(\Gamma)$ lies on a $g$-cycle of length at most $2d-2$  and let us define
$$N_{\Gamma}(b,i)=\{v\; | \; v\in V(\Gamma), d_{\Gamma}(b,v)=i\}, 0\leq i\leq d.$$
It is easy to see that $|N_{\Gamma}(b,0)|=1, |N_{\Gamma}(b,1)|=k$,  and $|N_{\Gamma}(b, i)|\leq k(k-1)^{i-1}$ for $2\leq i\leq d-2.$ Since $b$ lies on a $g$-cycle, where $g\leq 2d-2,$ we obtain $|N_{\Gamma}(b,d-1)|\leq k(k-1)^{d-2}-1.$ Hence $|N_{\Gamma}(b,d)|\leq k(k-1)^{d-1}-(k-1).$
This implies the inequality
 $$\delta(\Gamma) = M(k,d)-|V(\Gamma)| =(1+k+k(k-1)+\ldots +k(k-1)^{d-1})-(|N_{\Gamma}(b,0)|+|N_{\Gamma}(b,1)|+$$
 $$+\ldots+|N_{\Gamma}(b,d-1)|+|N_{\Gamma}(b,d)|)\geq k>c.$$

Since the above considered classes cover all graphs of maximum
degree $k$, and in each case we have been able to show that the
defect of the graphs in these classes is greater than $c$, it follows
that the defect of all graphs $\Gamma$ of degree at most $k$ and
even diameter $d > D $ is greater than $c$.

\qed

\begin{remark}

The above calculations allow us to estimate the defect
$\delta(\Gamma)$ of any $k$-regular graph $ \Gamma $
with sufficiently large even diameter $d$ and girth at least $2d-1$.
Recall that this girth assumption yields a lower bound on $|V(\Gamma)|$
which can be used to obtain an upper bound on the defect:
$\delta(\Gamma)=M(k,d)-|V(\Gamma)|\leq \frac{k(k-1)^{d}-2}{k-2}-\frac{k(k-1)^{d-1}-2}{k-2}=k(k-1)^{d-1}.$
Therefore, the defect $\delta(\Gamma)$ of any $k$-regular graph of
sufficiently large even diameter $d$ and of girth at least $2d-1$
belongs to one of the intervals
$[(d+1)\sqrt{k-1}^{d}-d\sqrt{k-1}^{d-1}, k(k-1)^{d-1}]$ or $[(d+1)\sqrt{k-1}^{d}+d\sqrt{k-1}^{d-1}, k(k-1)^{d-1}].$

This observation is related to
the concept of a generalized Moore graph \cite{CCMS}: A $k$-regular
graph $\Gamma$ of diameter $d$ and girth at least $2d-1$ is called
a {\em generalized Moore graph}. As argued above, the defect of
any generalized Moore graph is bounded from above by
$ k(k-1)^{d-1}$. For example, both Moore graphs and
bipartite Moore graphs are generalized Moore graphs, with the Moore
graphs having the defect $0$ and the bipartite Moore graphs having
the defect $ (k-1)^d $. It has been conjectured that the diameter of a generalized Moore graph cannot exceed $6$ (e.g., notes from a talk
delivered by L.K.\ J\o rgensen in Bandung, 2012). Our lower bounds on the defects of non-bipartite generalized Moore graphs do not seem
to contribute to the resolution of this conjecture.
\end{remark}

 \qed


\begin{thebibliography}{99}
{\small
\bibitem{AB} N. Alon, Eigenvalues and expanders, {\em Combinatorica} {\bf 6}
(1986) 83-96.
\bibitem{D} E.\ Bannai and T.\ Ito, On finite Moore graphs, {\em J. Fac. Sci. Tokyo Univ.} {\bf 20} (1973)
191-208.
\bibitem{B} J. -C.\ Bermond and B. Bollob\'as, The diameter of graphs: A survey, {\em Congressus
Numerantium} {\bf 32} (1981) 3-27.
\bibitem {L} N. I.\ Biggs, {\em Algebraic Graph Theory}, Cambridge University Press, Second Edition,
Great Britain (1993).
\bibitem{Big}
N. Biggs, Girth, valency, and excess,
{\em Linear Algebra Appl.} {\bf 31} (1980) 55-59.
\bibitem{CCMS}
V.G.\ Cerf, D.D.\ Cowan, R.C.\ Mullin, and R.G.\ Stanton,
Computer networks and generalized Moore graphs,
{\em Proc. 3rd Manitoba Conf. Numer. Math.}, Winnipeg 1973, (1974) 379-398.
\bibitem{C} R. M.\ Damerell, On Moore graphs, {\em Proc. Cambridge Phil. Soc.} 74 (1973) 227-236.
\bibitem{delorme} C. Delorme and G. P. Villavicencio, On graphs with cyclic defect or excess,
{\em Electron. J. Combin.} {\bf 17} (2010) no. 1.
\bibitem{A}  G.\ Exoo, R.\ Jajcay, M.\ Ma\v caj and J.\ \v Sir\' a\v n, On the defect of vertex-transitive graphs of given degree and diameter, submitted for publication.
\bibitem {F} W.\ Feit and G.\ Higman, The nonexistence of certain generalized polygons, {\em Journal
of Algebra} {\bf 1} (1964) 114-131.
\bibitem{Lov}
L. Lov\' asz, Eigenvalues of graphs, manuscript (2007).
\bibitem{ramanujan}  A. \ Lubotzky, R. \ Phillips and P.\ Sarnak, Ramanujan Graphs, {\em Combinatorica} {\bf 8} (3) (1988) 261-277.
\bibitem{K} M.\ Miller and J.\ \v Sir\' a\v n, Moore graphs and beyond: A survey of the degree/diameter problem, {\em Electron.  J. Combinatorics}, Dynamic Survey 14, 2005.
\bibitem {S} R. C.\ Singleton, On minimal graphs of maximum even girth, {\em J. of Comb. Theory} {\bf 1} (3) (1966) 306-332.
}

\end{thebibliography}
\end{document}